%% file: Version2.tex
\documentclass{amsart}
\title{A Vanishing Result for Toric Varieties Associated with Root Systems}

\usepackage{epic,eepic, graphicx, amssymb}

\newtheorem{definition}{Definition}[section]
\newtheorem{theorem}[definition]{Theorem}

\newtheorem{conjecture}[definition]{Conjecture}
\newtheorem{lemma}[definition]{Lemma}
\newtheorem{remark}[definition]{Remark}

\newcommand{\R}{\mathbb{R}}
\newcommand{\Z}{\mathbb{Z}}

\newcommand{\noadd}[1]{}

\input xy

\xyoption{all}

\newcommand{\Conv}{\operatorname{Conv}}

\author{Q\"endrim R. Gashi}
\begin{document}
\begin{abstract}
Consider a root system $R$ and the corresponding toric variety $V_R$ whose fan is the Weyl fan and whose lattice of characters is given by the root lattice for $R$. We prove the vanishing of the higher cohomology groups for certain line bundles on $V_R$ by proving a purely combinatorial result for root systems. These results are related to a converse to Mazur's Inequality for (simply-connected) split reductive groups.
\end{abstract}

\maketitle
\begin{center}

\end{center}

\section{Introduction}
 
The problem of $p$-adically estimating the number of points on an algebraic variety over a finite field of characteristic $p$ is quite old. An answer to a special case of this problem is given by the classical theorem of C. Chevalley and E. Warning, which asserts that if $G(x_1,...,x_n)$ is a polynomial of degree less than $n$ with integral coefficients, then the number of roots of $$G(x_1,...,x_n)\equiv 0 (\emph{mod p})$$ is divisible by $p$.

N. Katz conjectured a sharper $p$-adic estimate for the number of solutions to the above equation, which then B. Mazur proved in \cite{mazur} (and later P. Berthelot and A. Ogus completed in \cite{ogus}), making use of the then recently discovered crystalline cohomology. Mazur's result is now frequently referred to as Mazur's Inequality and it is most easily stated using the so-called Newton and Hodge vectors. We will recall this inequality in the setting of $F$-isocrystals.

An $F$-isocrystal is a pair $(N,F)$, where $N$ is a finite-dimensional vector space over the fraction field $K$ of the ring of Witt vectors $W(\overline{ \mathbb{F}}_p)$, equipped with a Frobenius-linear bijective endomorphism $F$ of $N$. 

Suppose now that our isocrystal $(N,F)$ is $n$-dimensional. If $M$ is a $W(\overline{ \mathbb{F}}_p)$-lattice in $N$, then we can associate to it the Hodge vector $\mu (M) \in \mathbb{Z}^n$, which measures the relative position of the lattices $M$ and $FM$. By Dieudonn\'{e}-Manin theory, we can associate to $N$ its Newton vector $\nu (N,F) \in \mathbb{Q}^n$, which classifies $F$-isocrystals of dimension $n$ up to isomorphism.

If $\geq$ stands for the ``usual dominance order'', then Mazur's Inequality asserts that $\mu (M) \geq \nu (N,F)$. In other words, if $\mu (M)=(\mu_1,...,\mu_n)$ and $\nu (M)=(\nu_1,...,\nu_n)$, then $\mu_1 \geq \nu_1$, $\mu_1 + \mu_2 \geq \nu_1+ \nu_2$,..., $\mu_1+...+\mu_{n-1} \geq \nu_1+...+\nu_{n-1}$, and $\mu_1+...+\mu_n = \nu_1+...+\nu_n$.

R. Kottwitz and M. Rapoport in \cite{krapo} proved a converse to this inequality. Namely, they proved that if we let $(N,F)$ be an isocrystal of dimension $n$, and let $\mu=(\mu_1,...,\mu_n) \in \mathbb{Z}^n$, with $\mu_1 \geq \mu_2 \geq ... \geq \mu_n$, be such that $\mu \geq \nu(N,F)$, then there exists a $W(\overline{ \mathbb{F}}_p)$-lattice $M$ in $N$ satisfying $\mu = \mu(M)$.

Both Mazur's Inequality and its converse can be regarded as statements for the group $GL_n$---the dominance order arises naturally in the context of the root system for $GL_n$. Kottwitz and Rapoport in \cite{krapo} (see also \cite{frapo} and \cite{rapritch}) formulated a group-theoretic version of the above converse to Mazur's Inequality, which they proved for $GL_n$ and $GSp_{2n}$. They also reduced this problem to a combinatorial-type problem formulated purely in terms of root systems (see e.g. \cite{kot} for an explicit statement). Then C. Lucarreli in \cite{cathy}, proved that combinatorial statement for all (split connected) classical groups.

In \cite{qendrim}, a new interpretation for the root-system combinatorial problem mentioned above was introduced; it was shown that it is equivalent to the vanishing of higher cohomology groups for certain line bundles on toric varieties associated with root systems. In loc.cit. a generalized version of this was proved for $GL_n$ and the usual version for $G_2$, in particular giving an easy proof for the converse to Mazur's Inequality for $GL_n$ and $G_2$ respectively. One of the surprising outcomes of these results, apart from the link of the converse to Mazur's Inequality to toric varieties, is that they improve the classical vanishing theorems for these toric varieties (compare, for example, results in \cite{mustata}). They also give new ideas for new vanishing theorems for general toric varieties. 

In the current paper, we take a new approach to proving the combinatorial problem mentioned above (see Theorem \ref{maintheorem} below). Also, unlike in other works related to this problem that have appeared so far, all the root systems are considered at once and the method of the proof applies to each of them. We prove a vanishing result (see Theorem \ref{torictheorem} below) for toric varieties associated with \emph{any} root system. Also, we believe that these results provide a crucial step and a clear strategy for proving the converse of Mazur's Inequality for all (split, connected) groups.

Finally, we mention that these toric varieties have appeared in the De Concini-Procesi theory of group compactifications (see e.g. \cite{brion}, pp.187--206) and are still actively used and studied in that field. They have also appeared on the work related to the local trace formula (see e.g. \cite{akot}, \S23) and the fundamental lemma (see e.g. \cite{claumon}, \S5).

\paragraph{\bf{Acknowledgments:}} The author thanks his thesis adviser, Professor Robert Kottwitz, for comments on an earlier draft of this paper, and for continuous encouragement and support.

\section{Set-up and Results}

We follow the terminology from \cite{bourbaki}. Let $R$ be an irreducible, reduced root system and let $Q(R)$ stand for the root lattice for $R$. Denote by $W_R$ the Weyl group for $R$. Let $x \in Q(R)$ and consider $O_x := \{ wx: w \in W_R \}$, the Weyl orbit of $x$; write $\Conv (O_x)$ for the convex hull of $O_x$ in $Q(R)\otimes_{\mathbb{Z}} \mathbb{R}$. Fix a root $\alpha$ in $R$ and denote by $\alpha^{\vee}$ the corresponding coroot. Suppose that $z=y+\frac{1}{2}\alpha \in Q(R)$ is such that $y\in \Conv (O_x)$ and $\langle y, \alpha^{\vee} \rangle = 0$.

The main result of this paper, formulated combinatorially, is the following:

\begin{theorem} With $x$ and $z$ as above, we have that $z\in\Conv\left( O_x \right).$
\label{maintheorem}
\end{theorem}

We prove this theorem in the next section but in the rest of this section we briefly explain how it is related to a converse to Mazur's inequality and to toric varieties associated with root systems (for more details on this relationship as well as for detailed references see e.g. \cite{qendrim}, Introduction and Section 1 in particular).

Once a notation is introduced, it will be fixed for the rest of the paper. For a short time we will not be working over complex numbers.

Let $F$ be a finite extension of $\mathbb{Q}_p$. Denote by $\mathcal{O}_F$ the ring of integers of $F$. Suppose $G$ is a split, connected reductive group, $B$ a Borel subgroup, and $T$ a maximal torus in $B$, all defined over $\mathcal{O}_F$. Let $P=MN$ be a parabolic subgroup of $G$ which contains $B$, where $M$ is the unique Levi subgroup of $P$ containing $T$.

We write $X$ for the set of cocharacters $X_*(T)$. Then $X_{G}$ and $X_{M}$ will stand for the quotient of $X$ by the coroot lattice for $G$ and $M$, respectively. Also, we let $\varphi_{G}:X\rightarrow X_G$ and $\varphi_{M}:X\rightarrow X_M$ denote the respective natural projection maps.

Let $\mu\in X$ be $G$-dominant and let $W$ be the Weyl group of $T$ in $G$. The group $W$ acts on $X$ and so we consider $W\mu :=\left\{w(\mu): w\in W\right\}$ and the convex hull of $W\mu$ in $\mathfrak{a} :=X\otimes_{\mathbb{Z}}\mathbb{R}$, which we denote by $\Conv\left(W\mu\right)$. Define
$$P_{\mu}=\left\{\nu\in X: (i)\, \varphi_{G}(\nu)=\varphi_{G}(\mu); \, \emph{and} \,\,(ii)\, \nu\in \Conv\left(W\mu\right)\right\}.$$

Let $\mathfrak{a}_M:=X_M\otimes_{\Z}\R$ and write $\emph{pr}_M:\mathfrak{a}\rightarrow\mathfrak{a}_M$ for the natural projection induced by $\varphi_M$. Note that $X_M$ is a quotient of $X$, but we can consider $\mathfrak{a}_M$ as a subspace of $\mathfrak{a}$ (after tensoring with $\mathbb{R}$ any possible torsion is lost).

We now recall an important conjecture of Kottwitz and Rapoport. See e.g. ~\cite{kot}, sections 4.3 and 4.4, where it is explained how a converse to Mazur's inequality follows from this conjecture.

\begin{conjecture} (Kottwitz-Rapoport) Let $G$ be a split, connected, reductive group over $F$. Keeping the same notation as above, we have

$$\varphi_M\left(P_\mu\right)=\left\{\nu_1\in X_M :(i)\,\nu_1, \mu \,\,\emph{have the same image in}\,\, X_G;\right.$$
$$\left. \hspace{5cm} (ii)\, \emph{the image of}\, \,\nu_1\, \emph{in} \,\mathfrak{a}_M \,\emph{lies in}\, pr_M\left(\Conv W\mu \right)\right\}.$$
\label{conjecture}
\end{conjecture}

It is worth mentioning that C. Lucarelli (see e.g. \cite{cathy}, Theorem 0.2) has proved that this conjecture holds for every split, connected reductive group over $F$ where every irreducible component of its Dynkin diagram is of type $A_n, B_n, C_n$, or $D_n$. In \cite{qendrim}, this has been further extended to include the group $G_2$ and a strengthened result has been proved in the case of the group $GL_n$ (see theorems A and B in loc. cit.). We again refer the reader to the Introduction in \cite{qendrim} where references are given for earlier work on this conjecture by other mathematicians.

In this set-up, our main result can be reformulated as follows.

\begin{theorem} The above conjecture of Kottwitz and Rapoport is true for all split, connected, simply-connected semisimple groups in the case when the parabolic subgroup $P$ is of semisimple rank 1.

\label{grouptheorem}
\end{theorem}

We would like to make a comment about some advantages and disadvantages of the proof of Theorem \ref{maintheorem} (and therefore Theorem \ref{grouptheorem}) . One of the advantages of the proof presented in this paper is that, unlike the previous proofs of (parts of) the Kottwitz-Rapoport Conjecture, all the root systems are considered at once and the method of the proof applies to each of them. Therefore, in particular, one gets a clearer picture, at least on the root-system level, of the entire problem of a converse to Mazur's inequality. Our approach has some serious disadvantages, however. The assumption that the group $G$ is simply-connected is used in an essential way (equivalently, we only consider the root lattice and not the weight lattice for our corresponding root system). Also, we only deal with the case when the parabolic subgroup $P$ has semisimple rank 1. Admittedly, the latter is less of a ``serious'' problem and, in fact, the author believes that it can be overcome without introducing new strategies in the proof.

Let us now  very briefly recall how this is all connected to vanishing results on certain toric varieties. The reader is encouraged to consult \cite{qendrim} (especially the Introduction and Section 1 therein) for more details. 

From now on we assume that $G$ is simply-connected. Then $\hat{G}$ is adjoint. Let $\hat{G}$ and $\hat{T}$ be the (Langlands) complex dual group for $G$ and $T$, respectively. Then the corresponding toric variety, which we denote by $V_G$, is given by requiring that its fan be given by the Weyl fan in $X_* (\hat{T}) \otimes_{\Z}\R$ and its torus be $\hat{T}$. (We would like to remark that the toric variety $V_G$ appears naturally in the theory of group compactifications---see e.g. ~\cite{brion}, pp.187--206.) Note that it is now safe for the reader to assume that we are working over complex numbers.

Using the canonical surjection, we define a map $$pr_M: X^*(\hat{T}) \twoheadrightarrow X^*(\hat{T})/R_{\hat{M}},$$ where $R_{\hat{M}}$ stands for the root lattice for $\hat{M}$ (note that the codomain of the last map is just $X^*(Z(\hat{M}))$).

For $M$ as above (a Levi subgroup containing $T$), we will need a toric variety $Y_M^G$ for the torus $Z(\hat{M})$ (see e.g. ~\cite{akot}, \S23.2). Note that $Z(\hat{M})$ is a subtorus of $\hat{T}$ and so $X_*(Z(\hat{M}))$ is a subgroup of $X_*(\hat{T})$. The collection of cones from the Weyl fan inside $X_*(\hat{T})\otimes_{\mathbb{Z}}\mathbb{R}$ that lie in the subspace $X_*(Z(\hat{M})) \otimes_{\mathbb{Z}}\mathbb{R}$ gives a fan. This is the fan for the nonsingular, projective toric variety $Y_M^G$.

Let us now assume that our parabolic subgroup $P$ is of semisimple rank 1. This implies that the root lattice $R_{\hat{M}}$ is just $\mathbb{Z} \alpha$ for a unique, up to a sign, root $\alpha$ of $\hat{G}$, and that the toric variety $Y_M^G$, which we now denote by $D_{\alpha}$, is a non-torus-invariant divisor in $V_G$. The map $pr_M$ will now be denoted by $p_{\alpha}$. By tensoring with $\mathbb{R}$ we get a map from $p_{\alpha}$, which we still denote by $$p_{\alpha}:X^*(\hat{T}) \otimes_{\mathbb{Z}}\mathbb{R} \twoheadrightarrow ( X^*(\hat{T})/\mathbb{Z}\alpha )\otimes_{\mathbb{Z}}\mathbb{R} .$$ Since tensoring with $\mathbb{R}$ will lose any possible torsion, we can identify the codomain of the last map with the coroot hyperplane $$[\alpha^{\vee} = 0] : =\{x \in  X^*(\hat{T})\otimes_{\mathbb{Z}}\mathbb{R}: \langle \alpha^{\vee},x \rangle =0  \} ,$$ where $\langle ,\rangle$ is the canonical pairing between cocharacters and characters, and $\alpha^{\vee}$ stands for the coroot of $\hat{G}$ corresponding to $\alpha$. Also, the fan of $D_{\alpha}$ is contained in the root hyperplane $$[\alpha = 0] : = \{ x \in X_*(\hat{T}) \otimes_{\mathbb{Z}}\mathbb{R}: \langle x, \alpha \rangle =0 \}.$$

Now let $\mathcal{L}$ be a $\hat{T}$-equivariant line bundle on $V_G$ that is generated by its sections. Also, assume that $\mathcal{L}$ is invariant under the obvious action of $W$. Then we have a short-exact sequence of sheaves on $V_G$: $$0\longrightarrow \mathcal{J}_{D_{\alpha}} \otimes \mathcal{L} \longrightarrow \mathcal{L} \longrightarrow i_*(\mathcal{L}|_{D_{\alpha}}) \longrightarrow 0,$$
where $\mathcal{J}_{D_{\alpha}}$ is the ideal sheaf of $D_{\alpha}$ and $i$ is the inclusion map $D_{\alpha} \hookrightarrow V_G$. Note that $H^i(V_G, \mathcal{L})$=0, for all $i>0$, since $\mathcal{L}$ is generated by its sections and $V_G$ is projective, and also $$H^i(V_G,i_*(\mathcal{L}|_{D_{\alpha}}))=H^i(D_{\alpha},\mathcal{L}|_{D_{\alpha}})=0,$$ for all $i>0$, since $\mathcal{L}|_{D_{\alpha}}$ is generated by its sections and $D_{\alpha}$ is a projective toric variety. Therefore the short-exact sequence gives rise to the long-exact sequence
$$... \longrightarrow H^0(V_G, \mathcal{L}) \stackrel{\varphi}{\longrightarrow} H^0(V_G,i_*(\mathcal{L}|_{D_{\alpha}})) \longrightarrow H^1(V_G,\mathcal{J}_{D_{\alpha}} \otimes \mathcal{L}) \longrightarrow 0.$$ Using the usual combinatorial interpretation of the 0-th cohomology (see e.g. \cite{fulton}, pg. 66) we see that the map $\varphi$ is induced by the map $p_{\alpha}$. Also, clearly, the surjectivity of $\varphi$ is equivalent to $H^1(V_G,\mathcal{J}_{D_{\alpha}} \otimes \mathcal{L})=0.$ Then one finds (as is explained in Section 1 of \cite{qendrim}) that $\varphi$ is surjective if and only if Theorem \ref{grouptheorem} is true. Thus, in this language, our result can be written as follows.

\begin{theorem}
With notation as above, we have that $$H^i(V_{G},\mathcal{J}_{D_{\alpha}} \otimes \mathcal{L})=0\, , \forall i >0.$$ 
\label{torictheorem}
\end{theorem}

Before we end this section it is worth mentioning that vanishing results like the one in Theorem \ref{torictheorem} are also of independent interest just from a toric-variety point of view. A very important vanishing result for toric varieties has been proved by Musta\c{t}\u{a} (see e.g. ~\cite{mustata}, Theorem \ref{maintheorem}), and for the toric varieties associated with root systems, Theorem \ref{torictheorem} gives the vanishing of higher cohomology groups for more line bundles on these varieties. 

As we mentioned earlier, the next section is devoted to proving Theorem \ref{maintheorem}.

\section{Proof of Theorem \ref{maintheorem}}

We recall some of the notation from the previous section. We follow the terminology and notation from \cite{bourbaki}. Let $R$ be an irreducible, reduced root system and let $Q(R)$ stand for the root lattice for $R$. Suppose $ \{ \alpha_i : i \in I \}$ is the set of simple roots (for some choice of a chamber for R) and let $s_{\alpha_i}$ be the simple reflection corresponding to $\alpha_i,\, i \in I$ (we will assume that $I:\,=\{ 1,\ldots, n\}$ for some natural number $n$); then $\alpha^{\vee}_i$, $i \in I$, stands for the coroot corresponding to $\alpha_i$. Moreover, let $\{ \varpi_{i}:  i \in I \}$ denote the set of fundamental coweights, where $\langle  \varpi_i , \alpha_i \rangle = 1$, and $\langle , \rangle$ is the standard pairing between coweights and weights for $R$. We write $W$ for the Weyl group of $R$ (note that earlier we were writing $W_R$ instead).

Let $x \in Q(R)$ and consider $O_x := \{ wx: w \in W \}$, the Weyl orbit of $x$; write $\Conv (O_x)$ for the convex hull of $O_x$ in $Q(R)\otimes_{\mathbb{Z}} \mathbb{R}$. Fix a root $\alpha$ in $R$ and denote by $\alpha^{\vee}$ the corresponding coroot. Suppose that $z=y+\frac{1}{2}\alpha \in Q(R)$ is such that $y\in \Conv (O_x)$ and $\langle y, \alpha^{\vee} \rangle = 0$.

We would like to prove that $z \in \Conv(O_x)$. But, before we start the proof, we make a few reductions. First, by possibly choosing a different chamber of $R$, we can assume that $\alpha$ is a simple root, and we will therefore denote the latter by $\alpha_{i_0}$, where $i_0$ is some element in $I$. Second, it is clear that we may take $x$ to be dominant, since we are considering the convex hull of its orbit $O_x$ under $W$, and for some $w\in W$ we must have that $wx \in O_x$ is dominant. Finally, we claim that it suffices to prove the theorem under the assumption that $y$ is dominant; indeed, if $y$ is not dominant, then we can find $w\in W$ such that $wy$ is dominant, and we can apply the theorem to $wy$ to find that $wz \in \Conv (O_x)$, and therefore $z \in \Conv (O_x)$. 

We can thus conclude that to prove Theorem \ref{maintheorem} is suffices to prove the following lemma.

\begin{lemma}
Let $x\in Q(R)$ be dominant. Fix a simple root $\alpha_{i_0} \in S$ and suppose that $z=y+\frac{1}{2}\alpha_{i_0} \in Q(R)$ is such that
\label{firstlemma}
\begin{itemize}

\item[(i)] $\langle y, \alpha_{i_0}^{\vee} \rangle = 0$
\item[(ii)] $y$ is dominant
\item[(iii)] $\langle y, \varpi_i \rangle \leq \langle x, \varpi_i \rangle$, $\forall i \in I$.

\end{itemize}
Then we have that $z\in\Conv\left( O_x \right).$

\end{lemma}

(Since $y$ is dominant, the condition (iii) in the lemma is equivalent to the assumption $y\in \Conv (O_x)$.)

The idea of the proof is as follows. First, notice that $$\Conv(O_x)=\{u\in Q(R)_{\mathbb{R}} \,| \langle wu,\varpi_i \rangle 
\leq \langle x,\varpi_i \rangle,\,\forall w\in W,\, \forall i\in I \}.$$ In fact, we know that there exists a unique $w_0 \in W$ such that $w_0z$ is dominant, 
and since $x$ is assumed to be dominant and $\Conv(O_x)$ is invariant under the action of $W$, we get that
\begin{equation}
z \in \Conv(O_x) \Longleftrightarrow \langle w_0z, \varpi_i \rangle \leq \langle x, \varpi_i \rangle, \,\forall i \in I.
\label{*}
\end{equation}

Thus, the idea of the proof will be to find an element $w_0 \in W$ so that $w_0z$ is dominant and then to prove that the inequalities on the right-hand side of (\ref{*}) are satisfied.

Now we begin with the proof itsef. From the way $z$ was defined, we see that $$\langle z,\varpi_i \rangle =\langle y,\varpi_i \rangle, \,\forall i\in I \setminus \{ i_0 \}, $$ and $$\langle z,\varpi_{i_0} \rangle =\langle y,\varpi_{i_0} \rangle + \frac{1}{2}.$$ Then condition (iii) immediately gives $$\langle z,\varpi_i \rangle \leq \langle x,\varpi_i \rangle, \,\forall i\in I \setminus \{ i_0 \}.$$ And, since $z\in Q(R)$, or more precisely since $\langle z,\varpi_{i_0} \rangle$ is an integer (this is where we are using the fact that we are working with the root lattice and not the weight lattice!), the relation (iii) also gives $\langle z,\varpi_{i_0} \rangle \leq \langle x,\varpi_{i_0} \rangle.$ Thus
\begin{equation}
\langle z,\varpi_i \rangle \leq \langle x,\varpi_i \rangle, \,\forall i\in I.
\label{zsmaller}
\end{equation}

In the case when $z$ is dominant (or equivalently when $w_0=1$), the last set of inequalities implies that the inequalities in (\ref{*}) are satisfied and so the assertion of our lemma is true.

But $z$ need not be dominant. We have $$\langle z,\alpha_{i}^{\vee} \rangle = \langle y,\alpha_{i}^{\vee} \rangle + \frac{1}{2} \langle
 \alpha_{i_0},\alpha_{i}^{\vee} \rangle, \, \forall i\in I.$$ and (i) implies that $\langle z,\alpha_{i_0}^{\vee} \rangle = 1$. However, $y$ being dominant, and having $\langle z,\alpha_{i}^{\vee} \rangle \in \mathbb{Z}$, we find that $z$ is not dominant if and only if

\begin{equation}
\exists \, i_1 \in I \setminus \{ i_0 \} \,\textup{s.t.} \, \langle \alpha_{i_0},\alpha_{i_1}^{\vee} \rangle <-1 \, \textup{and} \, \langle 
z,\alpha_{i_1}^{\vee} \rangle =-1. 
\label{i_1}
\end{equation}
Note that if $i_1$ as above exists, then it is unique. Clearly, for $i\in I \setminus \{ i_1 \}$, $\langle z, \alpha^{\vee}_i \rangle \geq 0.$

Since $\langle \alpha_{i_0},\alpha_{i_1}^{\vee} \rangle <-1$ can only happen for non-simply-laced root systems, we can conclude that $z$ is always dominant for a simply-laced root system R and, by what we wrote above, the Lemma is true for such an R.

More specifically, if we use the notation as in the diagrams below, (\ref{i_1}) holds only if:

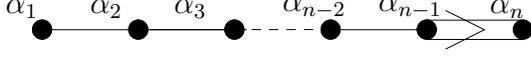
\input{Bn.eepic}

\begin{itemize}
\item[] $R=B_n$ and $\alpha_{i_0}=\alpha_{n-1}$, ${\alpha_{i_1}=\alpha_n}$, and $\langle z,\alpha_{i_1}^{\vee} \rangle =-1$; or
\end{itemize}

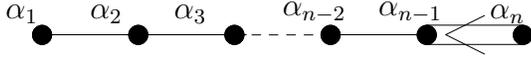
\input{Cn.eepic}

\begin{itemize}
\item[] $R=C_n$ and $\alpha_{i_0}=\alpha_n$, ${\alpha_{i_1}=\alpha_{n-1}}$, and $\langle z,\alpha_{i_1}^{\vee} \rangle =-1$; or
\end{itemize}

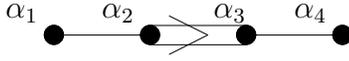
\input{F4.eepic}

\begin{itemize}
\item[] $R=F_4$ and $\alpha_{i_0}=\alpha_2$, ${\alpha_{i_1}=\alpha_3}$, and $\langle z,\alpha_{i_1}^{\vee} \rangle =-1$; or
\end{itemize}

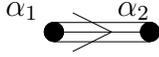
\input{G2.eepic}

\begin{itemize}
\item[] $R=G_2$ and $\alpha_{i_0}=\alpha_1$, ${\alpha_{i_1}=\alpha_2}$, and $\langle z,\alpha_{i_1}^{\vee} \rangle =-1$.
\end{itemize}

Let us now assume that (\ref{i_1}) holds. Our aim is to find $w_0 \in W$ so that $w_0z$ is dominant. We first apply the simple reflection 
$s_{\alpha_{i_1}}$ to bring $z$ ``closer'' to the dominant Weyl chamber. If we write $z_1:=s_{\alpha_{i_1}}(z),$ we see that 
$z_1=z+\alpha_{i_1}$, so $$\langle z_1,\alpha^{\vee}_{i}\rangle = \langle z,\alpha^{\vee}_{i} \rangle + \langle \alpha_{i_1},\alpha^{\vee}_{i} 
\rangle, \, \forall i\in I,$$
and it is also easy to check that $\langle z_1,\alpha^{\vee}_{i}\rangle \geq 0$, for $i=i_0,i_1$.

We then find that $z_1$ is not dominant if and only if
\begin{equation}
\exists \, i_2 \in I\setminus \{ i_0,i_1 \} \,\textup{s.t.} \, \langle \alpha_{i_1},\alpha_{i_2}^{\vee} \rangle =-1, \, \textup{and} \, 
\langle z,\alpha_{i_2}^{\vee} \rangle =0. 
\label{i_2}
\end{equation}
Let us remark that if $i_2$ as above exists, then it must be unique. Clearly, for $i \in I \setminus \{ i_2\},\,\langle z_1,\alpha^{\vee}_{i}\rangle \geq 0.$

Note that if $R=G_2$, then $z_1$ is always dominant, because if (\ref{i_2}) holds then $\alpha_{i_2}$ must be different from both $\alpha_{i_0}$ 
and $\alpha_{i_1}$, but this cannot happen since $G_2$ has only two simple roots. So, for $G_2$, $w_0=s_{\alpha_{i_1}}$.

For $R=B_n$, just as in the case of $G_2$, $z_1$ must be dominant, because there is no simple root $\alpha_{i_2}$, distinct from $\alpha_{i_0}$, 
such that $\langle \alpha_{i_1},\alpha_{i_2}^{\vee} \rangle =-1$. So, for $B_n$, $w_0=s_{\alpha_{i_1}}$.

If one uses the same notation as in the Dynkin diagrams above, then one can prove that the only possibility for $z_1$ to be non-dominant is if:
\begin{itemize}
\item[] $R=C_n$ and $\alpha_{i_0}=\alpha_n,\, {\alpha_{i_1}=\alpha_{n-1}}, \, {\alpha_{i_2}=\alpha_{n-2}}, \, \langle z,\alpha_{i_1}^{\vee} \rangle =-1,\, 
\langle z,\alpha_{i_2}^{\vee} \rangle =0$; or
\item[] $R=F_4$ and $\alpha_{i_0}=\alpha_2, \, {\alpha_{i_1}=\alpha_3},\, {\alpha_{i_2}=\alpha_4}, \, \langle z,\alpha_{i_1}^{\vee} \rangle =-1, \, \langle 
z,\alpha_{i_2}^{\vee} \rangle =0$.
\end{itemize}

Now we assume that (\ref{i_1}) and (\ref{i_2}) hold. We write $z_2:=s_{\alpha_{i_2}}(z_1),$ so $z_2=z+\alpha_{i_1}+\alpha_{i_2}$, and therefore 
$$\langle z_2,\alpha^{\vee}_{i}\rangle = \langle z,\alpha^{\vee}_{i} \rangle + \langle \alpha_{i_1},\alpha^{\vee}_{i} \rangle\ + 
\langle \alpha_{i_2},\alpha^{\vee}_{i} \rangle, \, \forall i\in I .$$
It is easy to see that $\langle z_2,\alpha^{\vee}_{i}\rangle \geq 0$, for $i=i_0,i_1,i_2$. We then get that $z_2$ is not dominant if and only if
\begin{equation}
\exists\, i_3 \in I\setminus \{ i_0,i_1,i_2\} \,\textup{s.t.} \, \langle \alpha_{i_2},\alpha_{i_3}^{\vee} \rangle =-1, \, \textup{and} \, \langle 
z,\alpha_{i_3}^{\vee} \rangle =0. 
\label{i_3}
\end{equation}
Note that if $i_3$ as above exists, then it is unique.  Clearly, for $i\in I\setminus \{ i_3\}, \, \langle z_2,\alpha^{\vee}_{i}\rangle \geq 0.$

Let us remark that this way we find that for $R=F_4$ we have $w_0=s_{\alpha_{i_2}} \circ s_{\alpha_{i_1}}$, because there is no $\alpha_{i_3}$ that satisfies (\ref{i_3}). 
So, it only remains to find a $w_0$ for $R=C_n$.

In any case, if $z_2$ is not dominant, we then apply $s_{\alpha_{i_3}}$ to it, to get $z_3:=s_{\alpha_{i_3}}(z_2)=z+\alpha_{i_1} + \alpha_{i_2} 
+\alpha_{i_3},$ and hence $$\langle z_3,\alpha^{\vee}_{i}\rangle = \langle z,\alpha^{\vee}_{i} \rangle + \langle \alpha_{i_1},\alpha^{\vee}_{i} \rangle
 + \langle \alpha_{i_2},\alpha^{\vee}_{i} \rangle + \langle \alpha_{i_3},\alpha^{\vee}_{i} \rangle, \, \forall i\in I .$$
Again $\langle z_3,\alpha^{\vee}_{i}\rangle \geq 0$, for $i=i_0,i_1,i_2, i_3$ and $z_3$ is not dominant if and only if
\begin{equation}
\exists \, i_4 \in I\setminus \{ i_0,i_1,i_2,i_3 \} \,\textup{s.t.} \, \langle \alpha_{i_3},\alpha_{i_4}^{\vee} \rangle =-1, \, \textup{and} \, \langle z,\alpha_{i_4}^{\vee} \rangle =0. 
\label{i_4}
\end{equation}
(If $i_4$ as above exists, then it must be unique.)

If (\ref{i_4}) holds, then we consider $z_4:=s_{\alpha_{i_4}}(z_3)=z+\alpha_{i_1} + \alpha_{i_2} +\alpha_{i_3}+ \alpha_{i_4}$. If $z_4$ is dominant, 
then the process of finding $w_0$ stops here, otherwise, we see that
\begin{equation}
\exists \, i_5 \in I\setminus \{ i_0,i_1,i_2,i_3,i_4 \} \, \textup{s.t.} \, \langle \alpha_{i_4},\alpha_{i_5}^{\vee} \rangle =-1, \, \textup{and} \, \langle z,\alpha_{i_5}^{\vee}
 \rangle =0.
\label{i_5}
\end{equation}
(If $i_5$ as above exists, then it must be unique.)

Since the number of simple roots is finite, we conclude, using induction, that $$\exists k\in \{1,\ldots , n-1 \}\, \textup{s.t.} \, w_0=s_{\alpha_{i_k}} \circ \ldots 
\circ s_{\alpha_{i_2}} \circ s_{\alpha_{i_1}},$$ and therefore $$w_0z=z+\alpha_{i_1}+\alpha_{i_2}+ \ldots+\alpha_{i_k}$$ is dominant. Clearly, by 
construction:
\begin{itemize}
\item[(a)] The $\alpha_{i_j}$'s appearing in $w_0$ are distinct;
\item[(b)] $\langle z,\alpha^{\vee}_{i_1} \rangle =-1$;
\item[(c)] $\langle z,\alpha^{\vee}_{i_j} \rangle =0$ for $j=2,\ldots,k$;
\item[(d)] $\langle \alpha_{i_j},\alpha^{\vee}_{i_{j+1}} \rangle =-1,$ for $j=1,\ldots,k-1$.

\end{itemize}

All that remains to finish the proof of Lemma \ref{firstlemma} (and therefore Theorem \ref{maintheorem}) is to check that $w_0z$ satisfies the right-hand side of (\ref{*}). This is done in the lemma below.

\begin{lemma} With notation as above \textup{(}in particular, $z\in Q(R)$ is assumed to satisfy \textup{(\ref{zsmaller})}\textup{)}, suppose that 
$w_0z=z+\alpha_{i_1}+\alpha_{i_2}+ \ldots+\alpha_{i_k}$ is dominant and conditions \textup{(a)}-\textup{(d)} hold. Then $$\langle 
w_0z,\varpi_i \rangle \leq \langle x,\varpi_i \rangle,\,\forall i\in I. $$
\end{lemma}

\begin{proof}
Conditions (\ref{zsmaller}) and (a) together imply that $$\langle w_0z,\varpi_i \rangle \leq \langle x,\varpi_i \rangle+1,\,\forall i\in \{i_1,\ldots,i_k \}, $$ 
and $$\langle w_0z,\varpi_i \rangle \leq \langle x,\varpi_i \rangle,\, \forall i\not \in \{i_1,\ldots,i_k \}.$$ Therefore, it suffices to check that $$\langle
 z,\varpi_i \rangle < \langle x,\varpi_i \rangle,\,\forall i\in \{i_1,\ldots,i_k \},$$ since $$\langle w_0z,\varpi_i \rangle = \langle z,\varpi_i \rangle + \sum^k_{j=1} \langle \alpha_{i_j},\varpi_i \rangle .$$

Suppose for a contradiction that $\langle z,\varpi_{i_1} \rangle = \langle x,\varpi_{i_1} \rangle$. Then, since $z \in Q(R)$ satisfies (\ref{zsmaller}), we 
get that $$z=x-\sum_{i\in I \setminus \{i_1\}}a_i\, \alpha_i,$$ for some non-negative integers $a_i$. But $x$ is dominant, so $\langle 
x,\alpha^{\vee}_{i_1} \rangle \geq 0$, and also $\langle \alpha_i , \alpha^{\vee}_{i_1} \rangle \leq 0,\,\forall i \neq i_1$. Thus, since $a_i$'s are 
non-negative, we have $\langle z, \alpha^{\vee}_{i_1} \rangle \geq 0$. But this contradicts our assumption (b) and therefore we must have $\langle 
z,\varpi_{i_1} \rangle < \langle x,\varpi_{i_1} \rangle$.

Assume now that $\langle z,\varpi_{i_m} \rangle = \langle x,\varpi_{i_m} \rangle$, for some $m\in \{2,\ldots,k \}$. Then, for similar reasons to those 
above, we can write $$z=x-\sum_{i\in I \setminus \{i_m\}}b_i\, \alpha_i,$$ for some non-negative integers $b_i$. Since $x$ is 
dominant, $\langle \alpha_i , \alpha^{\vee}_{i_m} \rangle \leq 0,\,\forall i \neq i_m$, and $b_i$'s are non-negative, we see that $b_i=0$, each 
time $\langle \alpha_i , \alpha^{\vee}_{i_m} \rangle < 0$, or else (c) would be contradicted. Condition (d) gives 
$\langle \alpha_{i_{m-1}},\alpha^{\vee}_{i_{m}} \rangle =-1$, so $b_{m-1}=0$.

We now get that $$z=x-\sum_{i\in I \setminus \{i_m,i_{m-1}\}}b_i\, \alpha_i,$$ for some non-negative integers $b_i$. From (d) we have 
$\langle \alpha_{i_{m-2}},\alpha^{\vee}_{i_{m-1}} \rangle =-1$, and so, because $x$ is dominant, $b_i$'s are non-negative, and the non-diagonal entries $\langle \alpha_i, \alpha^{\vee}_j \rangle \, (i \neq j)$ of the so-called Cartan matrix are not positive, we conclude, using (c), that $b_{i_{m-2}}=0$, or equivalently $$z=x-\sum_{i\in I \setminus \{i_m,i_{m-1}, i_{m-2}\}}b_i\, \alpha_i.$$

We continue this process by induction, to find that $$z=x-\sum_{i\in I \setminus \{i_m,i_{m-1},\ldots,i_1\}}b_i\, \alpha_i.$$ But this implies 
that $\langle z,\varpi_{i_1} \rangle = \langle x,\varpi_{i_1} \rangle$, contradicting the inequality $\langle z,\varpi_i \rangle < \langle x,\varpi_i \rangle $, 
demonstrated earlier.
\end{proof}

\begin{remark}
\emph{As was apparent in the proof of Lemma 1.1, $z$ is always dominant for simply-laced root systems, but it may fail to be dominant for the other root systems. However, when $z$ fails to be dominant, finding a $w_0\in W$ so that $w_0z$ is dominant was easier for some root systems than for others. More specifically, if $z$ is not dominant, then $w_0z=z+\alpha_{i_1}+\alpha_{i_2}+ \ldots+\alpha_{i_k}$ for certain distinct $\alpha_{i_j}$'s, and the number $k$ depends on the root systems. Call \emph{the defect of $z$} the number $k$ of simple reflections, as  in the proof of Lemma 1.1, required to make $z$ dominant. We then call \emph{the defect of the root system R}, denoted defect$(R)$, the maximum of defects of $z$, for $z\in Q(R)$. With this terminology, we have the following list, which, in a way, tells us the level of difficulty for solving our initial problem, stated in the Lemma 1.1.}

\begin{itemize}
\item[ 1.] defect $(A_n)=0$
\item[ 2.] defect $(D_n)=0$
\item[ 3.] defect $(E_6)=0$
\item[ 4.] defect $(E_7)=0$
\item[ 5.] defect $(E_8)=0$
\item[ 6.] defect $(B_n)=1$
\item[ 7.] defect $(G_2)=1$
\item[ 8.] defect $(F_4)=2$
\item[ 9.] defect $(C_n)=n-2$
\end{itemize}

\end{remark}

\end{document}

%% file: Bn.eepic
\setlength{\unitlength}{0.00083333in}
\begingroup\makeatletter\ifx\SetFigFont\undefined%
\gdef\SetFigFont#1#2#3#4#5{%
  \reset@font\fontsize{#1}{#2pt}%
  \fontfamily{#3}\fontseries{#4}\fontshape{#5}%
  \selectfont}%
\fi\endgroup%
{\renewcommand{\dashlinestretch}{30}
\begin{picture}(3293,378)(0,-10)
\put(225,132){\blacken\ellipse{120}{120}}
\put(225,132){\ellipse{120}{120}}
\put(825,132){\blacken\ellipse{120}{120}}
\put(825,132){\ellipse{120}{120}}
\put(1425,132){\blacken\ellipse{120}{120}}
\put(1425,132){\ellipse{120}{120}}
\put(2025,132){\blacken\ellipse{120}{120}}
\put(2025,132){\ellipse{120}{120}}
\put(2625,132){\blacken\ellipse{120}{120}}
\put(2625,132){\ellipse{120}{120}}
\put(3225,132){\blacken\ellipse{120}{120}}
\put(3225,132){\ellipse{120}{120}}
\path(900,132)(1500,132)
\path(825,132)(1425,132)
\path(225,132)(825,132)
\dashline{60.000}(1425,132)(2025,132)
\path(2025,132)(2625,132)
\path(2625,192)(3225,192)
\path(2625,72)(3225,72)
\path(2985,132)(2745,12)
\path(2745,252)(2985,132)
\put(525,207){\makebox(0,0)[lb]{{\SetFigFont{12}{14.4}{\rmdefault}{\mddefault}{\updefault}$\alpha_2$}}}
\put(1050,207){\makebox(0,0)[lb]{{\SetFigFont{12}{14.4}{\rmdefault}{\mddefault}{\updefault}$\alpha_3$}}}
\put(0,207){\makebox(0,0)[lb]{{\SetFigFont{12}{14.4}{\rmdefault}{\mddefault}{\updefault}$\alpha_1$}}}
\put(1725,207){\makebox(0,0)[lb]{{\SetFigFont{12}{14.4}{\rmdefault}{\mddefault}{\updefault}$\alpha_{n-2}$}}}
\put(2325,207){\makebox(0,0)[lb]{{\SetFigFont{12}{14.4}{\rmdefault}{\mddefault}{\updefault}$\alpha_{n-1}$}}}
\put(3021,207){\makebox(0,0)[lb]{{\SetFigFont{12}{14.4}{\rmdefault}{\mddefault}{\updefault}$\alpha_n$}}}
\end{picture}
}

%% file: Cn.eepic
\setlength{\unitlength}{0.00083333in}
\begingroup\makeatletter\ifx\SetFigFont\undefined%
\gdef\SetFigFont#1#2#3#4#5{%
  \reset@font\fontsize{#1}{#2pt}%
  \fontfamily{#3}\fontseries{#4}\fontshape{#5}%
  \selectfont}%
\fi\endgroup%
{\renewcommand{\dashlinestretch}{30}
\begin{picture}(3293,378)(0,-10)
\put(225,132){\blacken\ellipse{120}{120}}
\put(225,132){\ellipse{120}{120}}
\put(825,132){\blacken\ellipse{120}{120}}
\put(825,132){\ellipse{120}{120}}
\put(1425,132){\blacken\ellipse{120}{120}}
\put(1425,132){\ellipse{120}{120}}
\put(2025,132){\blacken\ellipse{120}{120}}
\put(2025,132){\ellipse{120}{120}}
\put(2625,132){\blacken\ellipse{120}{120}}
\put(2625,132){\ellipse{120}{120}}
\put(3225,132){\blacken\ellipse{120}{120}}
\put(3225,132){\ellipse{120}{120}}
\path(900,132)(1500,132)
\path(825,132)(1425,132)
\path(225,132)(825,132)
\dashline{60.000}(1425,132)(2025,132)
\path(2025,132)(2625,132)
\path(2625,192)(3225,192)
\path(2625,72)(3225,72)
\path(2745,132)(2985,252)
\path(2745,132)(2985,12)
\put(2325,207){\makebox(0,0)[lb]{{\SetFigFont{12}{14.4}{\rmdefault}{\mddefault}{\updefault}$\alpha_{n-1}$}}}
\put(525,207){\makebox(0,0)[lb]{{\SetFigFont{12}{14.4}{\rmdefault}{\mddefault}{\updefault}$\alpha_2$}}}
\put(1050,207){\makebox(0,0)[lb]{{\SetFigFont{12}{14.4}{\rmdefault}{\mddefault}{\updefault}$\alpha_3$}}}
\put(0,207){\makebox(0,0)[lb]{{\SetFigFont{12}{14.4}{\rmdefault}{\mddefault}{\updefault}$\alpha_1$}}}
\put(1725,207){\makebox(0,0)[lb]{{\SetFigFont{12}{14.4}{\rmdefault}{\mddefault}{\updefault}$\alpha_{n-2}$}}}
\put(3021,207){\makebox(0,0)[lb]{{\SetFigFont{12}{14.4}{\rmdefault}{\mddefault}{\updefault}$\alpha_n$}}}
\end{picture}
}

%% file: F4.eepic
\setlength{\unitlength}{0.00083333in}
\begingroup\makeatletter\ifx\SetFigFont\undefined%
\gdef\SetFigFont#1#2#3#4#5{%
  \reset@font\fontsize{#1}{#2pt}%
  \fontfamily{#3}\fontseries{#4}\fontshape{#5}%
  \selectfont}%
\fi\endgroup%
{\renewcommand{\dashlinestretch}{30}
\begin{picture}(2168,378)(0,-10)
\put(300,132){\blacken\ellipse{120}{120}}
\put(300,132){\ellipse{120}{120}}
\put(900,132){\blacken\ellipse{120}{120}}
\put(900,132){\ellipse{120}{120}}
\put(1500,132){\blacken\ellipse{120}{120}}
\put(1500,132){\ellipse{120}{120}}
\put(2100,132){\blacken\ellipse{120}{120}}
\put(2100,132){\ellipse{120}{120}}
\path(300,132)(900,132)
\path(900,192)(1500,192)
\path(900,72)(1500,72)
\path(1260,132)(1020,12)
\path(1020,252)(1260,132)
\path(1500,132)(2100,132)
\put(0,207){\makebox(0,0)[lb]{{\SetFigFont{12}{14.4}{\rmdefault}{\mddefault}{\updefault}$\alpha_1$}}}
\put(600,207){\makebox(0,0)[lb]{{\SetFigFont{12}{14.4}{\rmdefault}{\mddefault}{\updefault}$\alpha_2$}}}
\put(1800,207){\makebox(0,0)[lb]{{\SetFigFont{12}{14.4}{\rmdefault}{\mddefault}{\updefault}$\alpha_4$}}}
\put(1296,207){\makebox(0,0)[lb]{{\SetFigFont{12}{14.4}{\rmdefault}{\mddefault}{\updefault}$\alpha_3$}}}
\end{picture}
}

%% file: G2.eepic
\setlength{\unitlength}{0.00083333in}
\begingroup\makeatletter\ifx\SetFigFont\undefined%
\gdef\SetFigFont#1#2#3#4#5{%
  \reset@font\fontsize{#1}{#2pt}%
  \fontfamily{#3}\fontseries{#4}\fontshape{#5}%
  \selectfont}%
\fi\endgroup%
{\renewcommand{\dashlinestretch}{30}
\begin{picture}(968,378)(0,-10)
\put(300,132){\blacken\ellipse{120}{120}}
\put(300,132){\ellipse{120}{120}}
\put(900,132){\blacken\ellipse{120}{120}}
\put(900,132){\ellipse{120}{120}}
\path(300,192)(900,192)
\path(300,72)(900,72)
\path(660,132)(420,12)
\path(420,252)(660,132)
\path(300,132)(900,132)
\put(0,207){\makebox(0,0)[lb]{{\SetFigFont{12}{14.4}{\rmdefault}{\mddefault}{\updefault}$\alpha_1$}}}
\put(696,207){\makebox(0,0)[lb]{{\SetFigFont{12}{14.4}{\rmdefault}{\mddefault}{\updefault}$\alpha_2$}}}
\end{picture}
}